\documentclass[11pt,notitlepage,a4paper]{article}
\usepackage{graphicx}
\usepackage{amssymb}
\usepackage{amsmath}
\usepackage{appendix}

\def\tr{\mathop{\rm tr}\nolimits}

\def\tr{\mathop{\rm tr}\nolimits}

\newcommand{\bbeta}{\boldsymbol{\beta}}

\newcommand{\bdelta}{\boldsymbol{\delta}}
\newcommand{\bepsilon}{\boldsymbol{\epsilon}}

\newcommand{\bmu}{\boldsymbol{\mu}}

\newcommand{\bxi}{\boldsymbol{\xi}}

\newcommand{\bh}{\boldsymbol{h}}

\newcommand{\bu}{\boldsymbol{u}}

\newcommand{\bv}{\boldsymbol{v}}

\newcommand{\bW}{\boldsymbol{W}}
\newcommand{\bx}{\boldsymbol{x}}

\newcommand{\bY}{\boldsymbol{Y}}

\newcommand{\bzero}{\boldsymbol{0}}

\renewcommand{\L}{\mathcal{L}}

\def\tr{\mathop{\rm tr}\nolimits}

\def\tr{\mathop{\rm tr}\nolimits}

\newtheorem{theorem}{Theorem}

\newtheorem{lemma}[theorem]{Lemma}

\newtheorem{remark}{Remark}

\begin{document}

\title{\bf Restricted LASSO and Double Shrinking}

\bigskip

\author{{ M. Norouzirad$^1$, M. Arashi$^1$ and A.K.Md.Ehsanes Saleh$^2$}
\vspace{.5cm} \\\it $^{1}$ Department of Statistics, School of
Mathematical Sciences\\\vspace{.5cm} \it University of Shahrood, Shahrood, Iran\\
 \it $^2$ School of Mathematics and Statistics, Carleton University,\\\vspace{.5cm} \it Ottawa, Canada
 }

\date{}
\maketitle

\begin{quotation}
\noindent {\it Abstract:}
In the context of multiple regression model, suppose that the vector parameter of interest $\bbeta$ is subjected to lie in the subspace hypothesis $H\bbeta=\bh$, where this restriction is based on either additional information or prior knowledge. Then, the restricted estimator performs fairly well than the ordinary least squares one. In addition, when the number of variables is relatively large with respect to observations, the use of least absolute shrinkage and selection operator (LASSO) estimator is suggested for variable selection purposes.
In this paper, we define a restricted LASSO estimator and configure three classes of LASSO-type estimators to fulfill both variable selection and restricted estimation. Asymptotic performance of the proposed estimators are studied and a simulation is conducted to analyze asymptotic relative efficiencies. The application of our result is considered for the prostate dataset where the expected prediction errors and risks are compared. It has been shown that the proposed shrunken LASSO estimators, resulted from double shrinking methodology, perform better than the classical LASSO.
\par

\vspace{9pt} \noindent {\it Key words and phrases:} Double shrinking; LASSO; Preliminary test LASSO; Restricted LASSO; Stein-type Shrinkage LASSO.

\par

\vspace{9pt} \noindent {\it AMS Classification:} Primary: 62J07, 62F12,
Secondary:  62F30\par

\end{quotation}\par

\section{Introduction}
Consider a linear regression model with form
\begin{equation}
\bY_n= X_n\bbeta+\bepsilon_n, \label{eq:(2)}
\end{equation}
where the data are drawn  as $\{(\bx_i,Y_i)_{i=1}^n\}$, $\bx_i \in \mathbb{R}^p$ and $Y_i\in \mathbb{R}$ for $i=1,2,\ldots,n$. $X_{ij}$s are the regressors and $Y_i$ is the response variable of the $i$th observation and $\bbeta \in \mathbb{R}^p$ is unknown vector of coefficients to be estimated, $\bepsilon_n$ is the vector term with $E(\bepsilon_n)=0$ and $E(\bepsilon_n \bepsilon_n^T)=\sigma^2I_n(\sigma^2<\infty)$, $I_n$ is the identity matrix of order $n$.

In general, the main goal of the multiple regression model \eqref{eq:(2)} is the estimation of parameters and the prediction of response for a given design matrix. The estimation problem is usually solved through ordinary least squares (OLS) method where the parameters are estimated by the values minimizing the residual sum of squares $||\bY_n-X_n\bbeta||^2_2=\sum_{i=1}^{n} (Y_i-\bx_i^T\bbeta)^2$. Provided $X_n$ is of full rank, such that $X_n^TX_n$ is nonsingular and can be inverted, by the least squares method, the estimator of $\bbeta$ is written as
\begin{equation*}
\tilde{\bbeta}_n=(X_n^TX_n)^{-1}X_n^T\bY_n=C_n^{-1}X_n^T\bY_n;\qquad C_n=X_n^T X_n.
\end{equation*}
 The corresponding estimator of $\sigma^2$ is
\begin{equation}\label{eq:Se2}
s_e^2=\frac{1}{m}(\bY_n-X_n\tilde{\bbeta}_n)^T(\bY_n-X_n\tilde{\bbeta}_n); \quad m=n-p.
\end{equation}
It is obvious that $\tilde{\bbeta}_n\sim N_p(\bbeta,\sigma^2 C_n^{-1})$ independent of the distribution of $\frac{ms_e^2}{\sigma^2}$ which has a central chi-square distribution with $m$ d.f.

The standard procedures rely on the assumption that $C_n$ is nonsingular, otherwise $C_n$ cannot be inverted and the parameters cannot be uniquely estimated. The OLS regression method finds the unbiased linear combination of the $X_n$ that minimizes the residual sum of squares. However, if $p$ is large or the regression coefficients are highly correlated (multicollinearity), the OLS may yield estimates with large variance which reduces the accuracy of the prediction.  Possible solutions can be (1) \textit{variable selection}, for example best subset selection (Miller, 2002), (2) \textit{dimension reduction} techniques for example, principle component regression or partial linear regression; and (3) \textit{regularization} such that ridge (Hoerl and Kennard, 1970), LASSO (Tibshirani, 1996), SCAD (Zou and Hastie, 2005), elastic net (Fan and Li, 2001), and etc.

Variable selection is a method that  there are $p$ input variables, the objective is to select optimal model among all possible models. The most intuitive approach is maybe through preselction or subset selection. That is, to simply pick out a smaller subset of the covariates based on a certain relevant criterion and fit the (standard) model to these covariates only.  This approach (all possible models) is computationally infeasible, when $p$ is large (say, larger than $100$). There exist heuristics to cope with this problem such as forward selection, backward elimination or stepwise, but they are still \textit{unstable}. This procedure means that small changes of data result in large change of the estimator (Breiman, 1996). This method uses a hard decision rule (a variable survives or it dies).

The second approach was to use methods like principal components regression or partial least squares. These methods derive a small number of linear combinations of the original explanatory variables, and use these as covariates instead of the original variables. This may be reasonable for prediction purposes, but models are often difficult to interpret (Hastie et al., 2009).

Regularization methods are promising alternative. In these methods, coefficients are shrinkage rather than subset selection's result estimator. This process is more continues, and then we get lower variance than by subset selection and also reduce the prediction error of the full model. Shrinkage often improves prediction accuracy, trading off decreased variance for increased biased discussed in Hastie et al. (2009). These are also called shrinkage methods because they shrink the regression coefficients toward zero. The other name of this method is ``penalized regression methods" or more general, ``sparse regression".

A general sparse regression minimize the criterion
\begin{equation*}
f(\bbeta)+\sum_{j=1}^p P_\eta(|\beta_j|; \lambda),
\end{equation*}
where $f(\cdot)$ is a differentiable loss function (in linear regression usually $f(\bbeta)=||\bY_n-X_n\bbeta||_2^2$), $P(\cdot;\cdot)$ is a penalty function, $\lambda$ is a tuning parameter and $\eta$ index a penalty family. For example, we can refer to a power family, bridge regression (Frank and Friedman, 1993) as below:
\begin{equation*}
P_\eta(|\beta|,\lambda)=\lambda |\beta|^{\eta}; \qquad \eta\in[0,2].
\end{equation*}
If $\eta\in [0,1]$, then $P_\eta(|\beta|,\lambda)$ is concave, and it is a convex function if $\eta\in[1,2]$. Some special cases are (1) $\eta=0$, best subset regression, (2) $\eta=2$, ridge regression and (3) $\eta=1$, LASSO regression.

Elastic net family (Zou and Hastie, 2005) is
\begin{equation*}
P_\eta(|\beta|,\lambda)=\lambda\left\{(\eta-1)\frac{\beta^2}{2}+(2-\eta)|\beta|\right\}; \quad \eta\in [1,2].
\end{equation*}
In this family, $\eta=1$ result in the LASSO estimator and the ridge estimator is obtained by considering $\eta=2$.

Fan and Li (2001) defined the SCAD family as
\begin{equation*}
P_\eta(|\beta|,\lambda)=\left\{\begin{array}{ll}
\lambda|\beta|&|\beta|<\lambda\\
\lambda^2+\frac{\eta \lambda (|\beta|-\lambda)}{\eta-1}-\frac{\beta^2-\lambda^2}{2(\eta-1)}& |\beta|\in [\lambda,\eta \lambda]\\
\frac{\lambda^2(\eta+1)}{2}& |\beta|>\eta \lambda
\end{array}\right..
\end{equation*}
For small signals $|\beta|<\lambda$, it acts as LASSO, and for large signal $|\beta|>\eta\lambda$, the penalty flattens and lead to the unbiasedness of the regularized estimate.

Among all of the above regressions, the most famous is ridge. A disadvantage of this estimator is that the interpretation is not easy since the model includes all input variables. The LASSO is another important method. The $L_1$ penalty is used in LASSO while in ridge the $L_2$ is used. This tiny difference makes quantitative gaps practically as well as theoretically. The LASSO penalty shrinks each $\beta_j$ toward the origin and push irrelevant predictors to exact zero. Indeed, The LASSO can do variable selection and shrinkage estimation simultaneously. One very interesting property of LASSO is that the predictive model is sparse (i.e. some coefficient are exactly zero).

All regularization methods depend on one or more tuning parameters controlling the model complexity. Choosing the tanning parameters is an important part of the model fitting and is critical in statistical applications. There are two common used methods: (1) \textit{cross-validation}, and (2) \textit{information criteria} AIC and BIC offer good practice performance.
% that are given by
%\begin{eqnarray}
%AIC(\lambda)&=&\frac{||\bY-\hat{\bY}(\lambda)||_2^2}{\sigma^2}+2 df(\lambda),\cr
%\cr
%BIC(\lambda)&=&\frac{||\bY-\hat{\bY}(\lambda)||_2^2}{\sigma^2}+ \ln(n) df(\lambda),\nonumber
%\end{eqnarray}
%where $\hat{\bY}(\lambda)=X\hat{\bbeta}(\lambda)$ and $df(\lambda)$ is the effective degree of freedom of the selected model at $\lambda$.

%Efron (2004) used Stein (1981)'s theory of unbiased risk estimates and to show under differentiability condition on the mapping $\hat{\bY}(\lambda)$, we have
%\begin{equation*}
%df(\lambda)=\frac{1}{\sigma^2}\sum_{i=1}^{n} cov(\hat{Y_i}(\lambda),Y_i)=E\left[tr\left(\frac{\partial\hat{Y}(\lambda)}{\partial Y}\right)\right]
%\end{equation*}
%For least square estimate $df(\lambda)=tr[X(X^T X)^{-1}X^T]=p $, Zou et al. (2007) showed $df(\lambda)$ of the LASSO estimators is the number of non-zero coefficient.

In this paper, we focus only on the LASSO method. In the forthcoming section, we briefly introduce the LASSO estimator.

\subsection{LASSO Estimator}
Tibshirani (1996) proposed a new method for variable selection that produces an accurate, stable, and parsimonious model called LASSO (Least Absolute shrinkage and Selection Operator). The LASSO is a constrained version of OLS. Due to the sparseness property of the $L_1$ norm, the LASSO has been received much attention in recent years (Xu, 2014).  The ``LASSO" of Tibshirani (1996) is a least-squares problem regularized by $L_1$ norm, where we solve the following optimization problem
\begin{equation}\label{eq:LASSO}
\hat{\bbeta}^{L}=\text{min}_{\bbeta}\left\{\sum_{i=1}^n (y_i -\sum_{j=1}^p \beta_j X_{ij})^2\right\} \quad \text{subject to} \quad  \sum_{i=1}^p |\beta_j|\leq t ,
\end{equation}
where $t$ is a constant. If $t=0$, the model includes only the intercept term while the model becomes the full model when $t=\infty$.  If $t>\sum_{i=1}^p |{\beta}^*_j|$ where ${\bbeta}^*$ is the initial estimator of $\bbeta$ that usually it is considered as $\tilde{\bbeta}_n$, the OLS estimator, then the LASSO algorithm will yield the same as OLS estimate. However, if $0<t<\sum_{j=1}^p  |\tilde{\beta}_j|$, then the problem is equivalent to
\begin{equation}\label{eq:LASSO with CONSTRAINT}
\hat{\bbeta}^{L}=\arg \min_{\bbeta} \left\{\sum_{i=1}^n \left(y_i -\sum_{j=1}^p \beta_j X_{ij}\right)^2+\lambda_n \sum_{j=1}^p |\beta_j | \right\}, \qquad \lambda_n \ge 0,
\end{equation}
where $\lambda_n$ is a tunning parameter, controlling the level of sparsity in $\hat{\bbeta}^L$. The relation between $\lambda_n$ and LASSO parameter $t$ is one-to-one.

\subsection{Notation}
The following notation is used throughout the paper. For any dimension $d$, bold face letters denote vectors and normal face their elements, e.g. $\bv=(v_1,v_2,\ldots,v_d)^T$. Capital face letters denote matrices, e.g. $X$ and $\Sigma$. Let $A=(a_{ij})$ be an $m\times m$ matrix, then $A^T$ denoted the transpose of $A$, $\tr(A)=a_{11}+a_{22}+\ldots+a_{mm}$ is the trace of matrix $A$,  $A^{-1}=(a^{ij})$.

The $L_q$ norm of $v$ is $||\bv||_q=\left(\sum_{j=1}^d |v_j|^q\right)^{1/q}$ for $q>0$ and $|\cdot|$ to represent absolute value, applied. Design vectors, or columns of $X$, are denoted by $\bx_j$.

%For a real number $x$, let $sgn(x)$ be $-1$, $0$, or $1$ when $x<0$, $x=0$ and $x>0$; And for vector $\bv \in \mathbb{R}^p$, define $sgn(\bv)=(sgn(v_1),sgn(v_2),\ldots,sgn(v_p))^T$.

We will write $\text{min} \bv$ ($\text{max} \bv$) to denote the minimum (maximum) component of a vector $\bv$. Also, $\text{argmin} f(\cdot)$ ($\text{argmax} f(\cdot)$) is that component of the function's support which it result in minimum (maximum) value of $f$.

For notational convenience, we use $\Phi_p(\cdot;\bmu,\Sigma)$ and $\phi_p(\cdot;\bmu,\Sigma)$ to indicate the c.d.f. and the p.d.f. of the $p$-variate normal distribution with mean $\bmu$ and covariance matrix $\Sigma$, respectively. $H_\nu(\cdot;\Delta^2)$ denotes the c.d.f. of the $\chi^2$-distribution with degree of freedom $\nu$ and non centrality parameter $\Delta^2/2$.

Throughout, we use the following identities:
\begin{eqnarray}
&& E\left[\chi_{q+s}^{-2}(\Delta^2)\right]=E_r(q+s-2+2r)^{-1},\cr
&& E\left[\chi_{q+s}^{-4}(\Delta^2)\right]=E_r[(q+s-2+2r)(q+s-4+2r)]^{-1}\cr
&& E\left[\chi_{q+s}^{-2}(\Delta^2)I(\chi_{q+s}^2(\Delta^2)<k)\right]=E_r(q+s-2+2r)^{-1}H_{q+s-2+2r}(k;0),\cr
&& E\left[\chi_{q+s}^{-4}(\Delta^2)I(\chi_{q+s}^2(\Delta^2)<k)\right]=E_r[(q+s-2+2r)(q+s-4+2r)]^{-1} H_{q+s-4+2r}(k;0),\nonumber
\end{eqnarray}
where $E_r$ stands for the expectation with respect to a Poisson variable $r$ with parameter $\Delta^2/2$, and $I(A)$ is the indicator function of set $A$.

Finally, $\overset{\mathcal{P}}{\to}$ and $\overset{\mathcal{D}}{\to}$ are used to show convergence in probability and distribution, receptively.

We organize the paper as follows: In section \ref{sec:2}, a restricted LASSO estimator will be defined for inference under constraint and concept of double shrinking is introduced. Section \ref{sec:3} contains asymptotic distributions of the proposed estimators. In section \ref{sec:4}, a simulation study is conducted to analyze the relative efficiencies of the estimators, while an application of the results is considered for the well-known prostate dataset, where we compare expected prediction errors and asymptotic risk values.

\section{Restricted LASSO and Double shrinking}\label{sec:2}
‎Up to this point, it was assumed that the level of information had depend  on the sample, assuming no non-sample effect in estimation procedure. In this sense, we denote a LASSO estimator of $\bbeta$ by $\hat{\bbeta}_n^L$ and term it as unrestricted LASSO estimator (ULE).

‎However‎, in some situations it is possible to have some non-sample information (a priori restriction on the parameters) usually subjected to the model as constraints.

‎A set of $q$ linear restrictions on the vector $\bbeta$ can be written as $H\bbeta=\bh$‎. ‎Or we can suppose that our model is subjected to lie in the linear sub-space restriction
‎\begin{equation}\label{eq:Hb=h}‎
‎H \bbeta =\bh‎,
‎\end{equation}‎
where $H$ is a $q\times p$ $(q\leq p)$ matrix of known elements, with $q$ being the number of linear restriction to test, and $\bh$ is a $q\times 1$ vector of known components. ‎The rank of $H$ is $q$‎, which implies that the restrictions are linearly independent. ‎This restriction may be (a) a fact known from theoretical or experimental considerations, ‎(b) a hypothesis that may have to be tested or (c) an artificially imposed condition to reduce or eliminate redundancy in the description of the model (see Sengupta and Jammalamadaka‎, ‎2003)‎.

In this context, the LASSO estimator which satisfies \eqref{eq:Hb=h} will be  called  the restricted LASSO estimator (RLE), denoted by $\hat{\bbeta}_n^{RL}$. By the analogy of OLS estimator of $\bbeta$, subject to the restriction $H\bbeta=\bh$, we propose
\begin{equation}\label{eq:RLASSO}
\hat{\bbeta}_n^{RL}=\hat{\bbeta}_n^L-C_n^{-1}H^T(HC_n^{-1}H^T)^{-1}(H\hat{\bbeta}_n^L-\bh).
\end{equation}
When \eqref{eq:Hb=h} is satisfied, $\hat{\bbeta}_n^{RL}$ has smaller asymptotic risk than $\hat{\bbeta}_n^L$; However, for $H\bbeta\neq\bh$, $\hat{\bbeta}_n^{RL}$ may be biased and inconsistent in many cases. For this reason, it is plausible to follow Fisher's recipe and define a preliminary test LASSO estimator (PTLE) by taking $\hat{\bbeta}_n^L$ or $\hat{\bbeta}_n^{RL}$ according to acceptance or rejection of the null hypothesis
\[\mathcal{H}_o:H\bbeta=\bh. \]
This estimator will have the form
\begin{equation}\label{eq:PTLE}
\hat{\bbeta}_n^{PTL}=\hat{\bbeta}_n^L-(\hat{\bbeta}_n^L-\hat{\bbeta}_n^{RL})I(\mathcal{L}_n\leq \mathcal{L}_{n,\alpha}),
\end{equation}
where $\mathcal{L}_{n,\alpha}$ is the upper $\alpha$-level critical value of the exact distribution of the test   statistic $\mathcal{L}_n$ under $\mathcal{H}_o$. There will be two proposals for the test statistic $\mathcal{L}_n$. Following Saleh (2006) or Saleh et al. (2014), the test statistics is given by
\begin{equation}\label{eq:T.S.}
\mathcal{L}_n=\frac{(H\tilde{\bbeta}_n-\bh)^T(HC_n^{-1}H^T)(H\tilde{\bbeta}_n-\bh)}{s_e^2}.
\end{equation}
However, this test can be constructed upon the LASSO estimator. Here, we use the test in \eqref{eq:T.S.}. We believe that incorporating a test based on the LASSO estimator in analytical computations makes everything more easier.

The PTLE is highly dependent to the level of significance $\alpha$ and has discrete nature which simplifies to one of the extremes $\hat{\bbeta}_n^L$ or $\hat{\bbeta}_n^{RL}$ according to the output of the test. In this respect, making use of a continuous and $\alpha$-free estimator may make more sense. Now, we propose double shrinking idea which reflects a relevant estimator. It is well-known that the LASSO estimator shrinks coefficients toward the origin, however, when the restriction $H\bbeta=\bh$ is subjected to the model, it is of major importance that the estimator is shrunk toward the restricted one as well. Hence, there must be shrinking toward two directions or double shrinking, say. Consequently, we combine the idea of James-Stein (1961) shrinkage and LASSO to propose the following Stein-type shrinkage LASSO (SSLE) as
\begin{equation}\label{eq:SLE}
\hat{\bbeta}_n^{SSL}=\hat{\bbeta}_n^{L}-k_n(\hat{\bbeta}_n^L-\hat{\bbeta}_n^{RL})\mathcal{L}_n^{-1}, \qquad k_n=\frac{m(q-2)}{(m+2)},
\end{equation}
where $k_n$ is the shrinkage constant.

‎The estimator $\hat{\bbeta}^{SSL}_n$ may go past the estimator $\hat{\bbeta}_n^{RL}$‎. ‎Thus‎, we define the positive-rule Stein-type shrinkage LASSO estimator (PRSSLE) given by
\begin{eqnarray}‎
\hat{\bbeta}_n^{PRSSL}&=&\hat{\bbeta}^{RL}_n+\{1-k_n\mathcal{L}_n^{-1}\}I(\mathcal{L}_n>k_n)(\hat{\bbeta}^{L}_n-\hat{\bbeta}^{RL}_n),\cr‎
\cr
‎&=&\hat{\bbeta}^{SL}_n-(1-k_n\mathcal{L}_n^{-1})I(\mathcal{L}_n\leq k_n)(\hat{\bbeta}^{L}_n-\hat{\bbeta}^{RL}_n).\label{eq:PRSLE}‎
‎\end{eqnarray}
We note that, as the test based on $\mathcal{L}_n$ is consistent against fixed $\bbeta$ such that $H\bbeta\neq \bh$, the PTLE, SSLE and PRSSLE are asymptotically equivalent to the ULE for fixed alternative. Hence, we will investigate the asymptotic risks under local alternatives and compare the respective performance of the estimators.

\section{Asymptotic Distribution of the Estimators}\label{sec:3}
In sequel, the following regularity assumptions will be needed.
\begin{description}
\item[A1:] $\max_{1\leq i\leq n} \bx_i^T C_n^{-1}\bx_i\to 0$ as $n\to \infty$ where  $\bx_i^T$ is the $i$th row of design matrix $X$.
\item[A2:] $\lim_{n\to \infty} n^{-1}C_n=C$, where $C$ is finite and positive-definite matrix.
\end{description}
And also we need to notice the restriction $H\bbeta=\bh$ is not exact; rather, it is of the form $H\bbeta=\bh+\bxi$. In many situations, the asymptotic distribution of $\sqrt{n}s_e^{-1}(\bbeta_n^*-\bbeta)$ is equivalent to the $\sqrt{n}\sigma^{-1}(\tilde{\bbeta}_n-\bbeta)$ distribution as $n\to \infty$ under fixed alternatives,
\begin{equation*}
K_{\bxi}: H\bbeta=\bh+\bxi,
\end{equation*}
where $ \bbeta^*$ is an estimator of $\bbeta$. Then, to obtain the asymptotic distribution of $\sqrt{n}s_e^2(\bbeta_n^*-\bbeta)$, we consider the class of local alternatives, ${K_{(n)}}$ defined by
\begin{equation*}
K_{(n)}: H\bbeta=\bh+ n^{-\frac12}\bxi.
\end{equation*}
Now, let the asymptotic cumulative distribution function (c.d.f.) of $\sqrt{n}s_e^{-1}(\bbeta_n^*-\bbeta)$ under ${K_{(n)}}$ be
\begin{equation*}
G_p(\bx)=\lim_{n\to\infty}P_{K_{(n)}}\left\{\sqrt{n}s_e^{-1}(\bbeta_n^*-\bbeta)\leq \bx \right\}
\end{equation*}
If the asymptotic c.d.f. exists, then the asymptotic distributional bias (ADB) and quadratic bias (ADQB) are given by
\begin{equation*}
b(\bbeta^*)=\lim_{n\to\infty} E\left[\sqrt{n}(\bbeta_n^*-\bbeta)\right]=\int \bx dG_p(\bx)
\end{equation*}
and
\begin{equation*}
B(\bbeta_n^*)=\sigma^{-2}[b(\bbeta^*)]^T C [b(\bbeta_n^*)]
\end{equation*}
respectively, where $\sigma^2C^{-1}$ is the MSE-matrix of $\tilde{\bbeta}_n$ as $n\to\infty$. Defining
\begin{equation*}
M(\bbeta_n^*)=\int \bx \bx^T dG_p(x) = \lim_{n\to \infty} E\left[n(\bbeta_n^*-\bbeta)(\bbeta_n^*-\bbeta)^T\right],
\end{equation*}
as the asymptotic distributional MSE (ADMSE), we have the weighted risk of $\bbeta_n^*$ given by
\begin{equation*}
R(\bbeta_n^*)=\tr[M(\bbeta_n^*)]=\lim_{n\to\infty} E[n(\bbeta_n^*-\bbeta)^T(\bbeta_n^*-\bbeta)]
\end{equation*}
as the asymptotic distributional quadratic risk (ADQR).

For the proof of all following results, refer to the Appendix.
\begin{theorem}\label{thm:4}
Under $K_{(n)}: H\bbeta=\bh+n^{-\frac12}\bxi$ and the regularity assumptions, we have the following as $n\to \infty$,
\begin{description}
\item[(i)]
If $C$ is a nonsingular matrix and $\lambda_n/n\to\lambda_0\geq0$, then $\hat{\bbeta}_n^{L}\overset{\mathcal{P}}{\to}\text{argmin}(Z)$  where
\begin{equation*}
Z(\boldsymbol{\phi})=(\boldsymbol{\phi}-\bbeta)^T C (\boldsymbol{\phi}-\bbeta)+\lambda_0 \sum_{j=1}^{p}|\phi_j|.
\end{equation*}
\item[(ii)]$\hat{\bbeta}^{RL}_n\overset{\mathcal{P}}{\to}\text{argmin}(Z)-C^{-1}H^T(H C^{-1}H^T)^{-1}(H\text{argmin}(Z)-\bh)$.
\item[(iii)]
$\hat{\bbeta}^{L}_n-\hat{\bbeta}^{RL}_n\overset{\mathcal{P}}{\to} C^{-1}H^T(H C^{-1}H^T)^{-1}(H\text{argmin}(Z)-\bh)$.
\item[(iv)]
$\hat{\bbeta}_n^{PTL}\overset{\mathcal{P}}{\to} \text{argmin}(Z)-C^{-1}H^T(H C^{-1}H^T)^{-1}(H\text{argmin}(Z)-\bh) I(\L<\L_{\alpha})$
where $\L=\sigma^{-2}(H W+\bxi)^T(H C^{-1}H^T)^{-1} (HW+\bxi)$ and $\L_{\alpha}$ is the upper critical value of chi-squared distribution with $q$ d.f.
\item[(iv)]
$\hat{\bbeta}_n^{SSL}\overset{\mathcal{D}}{\to} \text{argmin}(Z)-kC^{-1}H^T(H C^{-1}H^T)^{-1}(H\text{argmin}(Z)-\bh)\L^{-1}$
\item[(v)]{\small
$\hat{\bbeta}_n^{PRSSL} \overset{\mathcal{D}}{\to} \text{argmin}(Z)-\left\{k\L^{-1}+(1-k\L^{-1})I(\L<k)\right\}C^{-1}H^T(H C^{-1}H^T)^{-1}(H\text{argmin}(Z)-\bh)$}
\end{description}
\end{theorem}

\begin{theorem}\label{thm:6}
Under the class of local alternatives, $\{K_{(n)}\}$, and regularity assumptions, we have the following as $n\to \infty$,
\begin{description}
\item[(i)] $\sqrt{n}(\tilde{\bbeta}_n-\bbeta)\sim N_p(\bzero,\sigma^2C^{-1})$.
\item[(ii)]
If $\frac{\lambda_n}{\sqrt{n}}\to\lambda_0\geq0$ and $C$ is a nonsingular matrix, then
\[ \sqrt{n}(\hat{\bbeta}^{L}_n-\bbeta)\overset{\mathcal{D}}{\to} \text{argmin}(V) \]
where
\begin{equation*}\label{eq:V(u)}
V(u)=-2\bu^T\bW+\bu^TC\bu+\lambda_0 \sum_{j=1}^{p}[u_jsgn(\beta_j)I(\beta_j\neq0)+|u_j|I(\beta_j=0)],
\end{equation*}
and $\bW\sim N_p(\boldsymbol{0},\sigma^2C)$.
\item[(iii)]
$\sqrt{n}(\hat{\bbeta}^{RL}-\bbeta)\overset{\mathcal{D}}{\to} \text{argmin}(V)-C^{-1}H^{T}(H C^{-1}H^T)^{-1}(H\text{argmin}(V)+\bxi)$.
\item[(iv)]
$\sqrt{n}(\hat{\bbeta}^{L}-\hat{\bbeta}^{RL}_n)\overset{\mathcal{D}}{\to} C^{-1}H^T(H C^{-1}H^T)^{-1}(H \text{argmin}(V)+\bxi)$
\item[(v)]
$\lim_{n\to \infty} P(\L_n\leq x)=H_q(x;\Delta^2)$ where $H_q(\cdot;\Delta^2)$ is the c.d.f. of non central chi squared distribution.
\item[(vi)]
$\sqrt{n}(\hat{\bbeta}^{PTL}-\bbeta)\overset{\mathcal{D}}{\to} \text{argmin(V)}-C^{-1}H^T(H C^{-1}H^T)^{-1}(H \text{argmin} (V)+\bxi)I(\L \leq \L_{\alpha})$.
\item[(vii)]
$\sqrt{n}(\hat{\bbeta}^{SSL}-\bbeta)\overset{\mathcal{D}}{\to} \text{argmin(V)}-k C^{-1}H^T(H C^{-1}H^T)^{-1}(H \text{argmin}(V)+\bxi)\L^{-1},\  k=(q-2)$.
\item[(viii)]
$\sqrt{n}(\hat{\bbeta}^{PRSSL}-\bbeta)\overset{\mathcal{D}}{\to}\text{argmin}(V)-\left\{k\L^{-1}+(1-k\L ^{-1})I(\L<k)\right\}C^{-1}H^{T}(H C^{-1}H^T)^{-1}$\\$\times (H\text{argmin}(V)+\bxi)$
\end{description}
\end{theorem}

\subsection{Null-Consistent Estimators}
In this section, suppose the LASSO is weakly consistent, i.e., $\lambda_n=o(n^{\tfrac12})$. Up to this point, we implemented a test statistic based on the OLS estimator, however, constructing a test based on the LASSO estimator will give the same asymptotic behavior in our setup. A test statistic based on the LASSO estimator will have form
\begin{equation}\label{eq:test-statistic}
\mathcal{L}_n=\frac{(H\hat{\bbeta}_n^L-\bbeta)^T(HC^{-1}H^T)^{-1}(H\hat{\bbeta}_n^L-\bbeta)}{s_L^2},
\end{equation}
where
\begin{equation}\label{eq:var-LASSO}
s_L^2=\frac{1}{m}(Y-X\hat{\bbeta}_n^L)^T(Y-X\hat{\bbeta}_n^L), \qquad m=n-p.
\end{equation}

\begin{theorem}\label{thm:3}
Under regularity assumptions and also in the class of local alternatives ${K_{(n)}}$,  $\L_n$, the likelihood ratio test statistics, converges in distribution to $\L$ which has the non central chi square distribution with $q$ d.f. and non centrally parameter $\Delta^2=\sigma^{-2}\bxi^T(HC^{-1}\bxi)^{-1}\bxi=\sigma^{-2}\bdelta^T C \bdelta$ where $\bdelta=C^{-1}H^T(H C^{-1} H^T)^{-1}\bxi$  and it defined as
\begin{equation*}
\L=\frac{(HW+\bxi)^T(HC^{-1}H^T)^{-1}(HW+\bxi)}{\sigma^2}
\end{equation*}
Where $W\sim N(\bzero,\sigma^2 C)$.
\end{theorem}

\begin{theorem}\label{thm:5}
In theorem \ref{thm:4},  if $\lambda_n=o(n)$, we have the following results,
\begin{description}
\item[(i)] $\text{argmin}(Z)=\bbeta$ and so $\hat{\bbeta}_n^{L}$ is consistent.
\item[(ii)] $\hat{\bbeta}^{RL}_n\overset{\mathcal{P}}{\to}\bbeta-\bdelta; \quad \bdelta=C^{-1}H^T(H C^{-1}H^T)^{-1}(H\bbeta-\bh)$.
\item[(iii)] $\hat{\bbeta}_n^{PTL}\overset{\mathcal{P}}{\to} \bbeta-\bdelta I(\L<\L_{\alpha})$.
\item[(iv)] $\hat{\bbeta}_n^{SL}\overset{\mathcal{P}}{\to}\bbeta-\bdelta\L^{-1}$.
\item[(v)]
$\hat{\bbeta}_n^{PRSL}\overset{\mathcal{P}}{\to}\bbeta-\left\{k\L^{-1}+(1-k\L^{-1})I(\L<k)\right\}\bdelta$.
\end{description}
\end{theorem}
\begin{remark}
According to Theorem \ref{thm:4}, under $H_0$, all estimators are consistent for $\bbeta$.
\end{remark}
In all the following results, proofs are directly deduced using the utilities in Saleh (2006) after some algebra.
\begin{theorem}\label{thm:7}
In Theorem \ref{thm:6},  if $\lambda_n=o(n^{\tfrac12})$, we have the following results:
\begin{description}
\item[(i)]   $W=\sqrt{n}(\tilde{\bbeta}_n-\bbeta)\overset{\mathcal{D}}{\to}  N_p(\boldsymbol{0},\sigma^2 C^{-1})$.
\item[(ii)]   $W_n^{(1)}=\sqrt{n}(\hat{\bbeta}_n^{LE}-\bbeta)\overset{\mathcal{D}}{\to}  N_p(\boldsymbol{0},\sigma^2 C^{-1})$. i.e. $W_n^{(1)}\overset{\mathcal{D}}{=} W$.
\item[(iii)] $W_n^{(2)}=\sqrt{n}(\hat{\bbeta}^{RL}_n-\bbeta)\overset{\mathcal{D}}{\to} N_p(-\bdelta,\sigma^2 A)$ where $\bdelta=C^{-1}H^T(H C^{-1}H^T)^{-1}\bxi$ and\\ $A=C^{-1}-C^{-1}H^T(H C^{-1}H^T)^{-1}H C^{-1}$.
\item[(iv)] $W_n^{(3)}=\sqrt{n}(\hat{\bbeta}^{L}_n-\hat{\bbeta}^{RL}_n)\overset{\mathcal{D}}{\to} N_p(\bdelta,\sigma^2 (C^{-1}-A))$ .
\item[(v)]$W_n^{(4)}=H\hat{\bbeta}^{L}_n-\bh\overset{\mathcal{D}}{\to} N_q(H\bbeta-\bh,\sigma^2 (HC^{-1}H^T))$,
\item[(vi)]
$\left[\begin{matrix}
W_n^{(1)}\\
W_n^{(3)}
\end{matrix}\right]\overset{\mathcal{D}}{\to} N_{2p}\left(\left[\begin{matrix}
\bzero\\
\bdelta
\end{matrix}\right],\sigma^2\left[\begin{matrix}
C^{-1}&C^{-1}-A\\
C^{-1}-A & C^{-1}-A
\end{matrix}\right]\right)$
\item[(vii)]$\left[\begin{matrix}
W_n^{(2)}\\
W_n^{(3)}
\end{matrix}\right]\overset{\mathcal{D}}{\to} N_{2p}\left(\left[\begin{matrix}
\bdelta\\
-\bdelta
\end{matrix}\right],\sigma^2\left[\begin{matrix}
A & 0\\
0 & C^{-1}-A
\end{matrix}\right]\right)$
\item[(viii)]$\left[\begin{matrix}
W_n^{(1)}\\
W_n^{(4)}
\end{matrix}\right]\overset{\mathcal{D}}{\to} N_{p+q}\left(\left[\begin{matrix}
0\\
H\bbeta-\bh
\end{matrix}\right],\sigma^2\left[\begin{matrix}
C^{-1} & C^{-1}H^T\\
H C^{-1} & HC^{-1}H^T
\end{matrix}\right]\right)$
%\item[(viii)] $\lim_{n\to \infty}P\{\sqrt{n}(\hat{\bbeta}_n^{PTL}-\bbeta)\leq \bx\}= H_q(\chi_q^2(\alpha);\Delta^2)\Phi_p(\bx+\bdelta,\bzero,\sigma^2 \bA)+(1-H_q(\chi_q^2(\alpha);\Delta^2))$\\
%$\times\Phi_p(\bx,\bzero,\sigma^2 C^{-1})$
%  where $\Phi_p(\cdot;\bmu,\bSigma)$ is the cdf of a $p$-variate normal distribution with mean $\bmu$ and covariance matrix, $\bSigma$, and $H_\nu(\cdot;\Delta^2)$ is the cdf of a noncentral chi-squared distribution with $\nu$ d.f. and noncentrality parameter $\Delta^2/2$.
\item[(ix)]$\sqrt{n}(\hat{\bbeta}^{SSL}-\bbeta)\overset{\mathcal{D}}{=} W-k\left\{\frac{C^{-1}H^T(HC^{-1}H^T)^{-1}(HW+\bxi)}{\sigma^{-2}(HW+\bxi)^T(HC^{-1}H^T)^{-1}(HW+\bxi)}\right\}$.
\item[(x)]
\begin{eqnarray}
\sqrt{n}(\hat{\bbeta}^{PRSL}-\bbeta)&\overset{\mathcal{D}}{=}& W-k\left\{\frac{C^{-1}H^T(HC^{-1}H^T)^{-1}(HW+\bxi)}{\sigma^{-2}(HW+\bxi)^T(HC^{-1}H^T)^{-1}(HW+\bxi)}\right\}\cr && +C^{-1}H^T(HC^{-1}H^T)^{-1}(HW+\bxi)\cr
&& \times \left\{1-\frac{k}{\sigma^{-2}(HW+\bxi)^T(HC^{-1}H^T)^{-1}(HW+\bxi)}\right\}I(\L<k).\nonumber
\end{eqnarray}
\end{description}
\end{theorem}

\begin{lemma} (Saleh, 2006)
We have the following identities:
\begin{itemize}
\item[(i)] $E\left[\chi_{q+2}^{-2}(\Delta^2)\right]=\exp{-\frac{\Delta^2}{2}}\sum_{r\geq0}\frac{1}{r!}\left(\frac{\Delta^2}{2}\right)^2\frac{1}{q+2r}=E_r[(q+2r)^{-1}]$ where $E_r$ stands for the expression with respect to the Poisson variable $r$ with mean $\frac{\Delta^2}{2}$.
\item[(ii)]$E\left[\chi_{q+2}^{-4}(\Delta^2)\right]=\exp{-\frac{\Delta^2}{2}}\sum_{r\geq0}\frac{1}{r!}\left(\frac{\Delta^2}{2}\right)^2\frac{1}{(q+2r)(q-2+2r)}=E_r[(q+2r)(q-2+2r)]^{-1}$.
\item[(iii)]$E\left[\chi_{q+2}^{-2}(\Delta^2)I(\chi_a^2(\Delta^2)<c)\right]=E_r[(q+2r)^{-1}] H_{q+2r}(c;0)$,
\item[(iv)]$E\left[\chi_{q+2}^{-4}(\Delta^2)I(\chi_a^2(\Delta^2)<c)\right]=E_r[(q+2r)^{-1}] H_{q+2r}(c;0)$.
\end{itemize}
\end{lemma}

\begin{theorem}\label{thm:9}
Suppose that all estimators are $\sqrt{n}$-consistent. Then, the ADB, ADQB, ADMSE and ADQR of the estimators are given by
\begin{eqnarray}
(i)&& b_1(\hat{\bbeta}_n^L)=\bzero,\qquad B_1(\hat{\bbeta}_n^L)=0, \qquad R_1(\hat{\bbeta}_n^L;W)=\sigma^2 \tr(WC^{-1}),\,\text{and}\,\quad  M_1(\hat{\bbeta}_n^L)=\sigma^2 C^{-1},\cr
\cr
(ii) && b_2(\hat{\bbeta}_n^{RL})=-\bdelta, \qquad B_2(\hat{\bbeta}_n^{RL})=\Delta^2, \qquad R_2(\hat{\bbeta}_n^{RL};W)=\sigma^2 \tr(W(C^{-1}-A))+\bdelta^T W \bdelta,\cr
&&\text{and}\qquad M_2(\hat{\bbeta}_n^{RL})=\sigma^2 (C^{-1}-A)+\bdelta\bdelta^T,\cr
\cr
(iii)&& b_3(\hat{\bbeta}_n^{PTL})=-\bdelta H_{q+2}\left(\chi_{q}^2(\alpha);\Delta^2\right), \qquad  B_3(\hat{\bbeta}_n^{PTL})=\Delta^2 \{H_{q+2}(\chi_q^2(\alpha);\Delta^2)\}^2,\cr
&& R_3(\hat{\bbeta}_n^{PTL};W)=\sigma^2 \tr(WC^{-1})-\sigma^2 \tr(W(C^{-1}-A))H_{q+2}(\chi_q^2(\alpha);\Delta^2)+\bdelta^T W\bdelta Z(\alpha;\Delta^2),\cr
&& M_3(\hat{\bbeta}_n^{PTL})=\sigma^2 C^{-1} -\sigma^2 (C^{-1}-A)H_{q+2}(\chi_q^2(\alpha);\Delta^2)+\bdelta \bdelta^T Z(\alpha;\Delta^2),\cr
\cr
(iv)&& b_4(\hat{\bbeta}_n^{SL})=-k \bdelta E\left[\chi_{q+2}^{-2}(\Delta^2)\right], \quad \text{where}\quad  k=\lim\limits_{n\to\infty}k_n=q-2,\cr
&&  B_4(\hat{\bbeta}_n^{SL})=k^2\Delta^2 \left\{E\left[\chi_{q+2}^{-2}(\Delta^2)\right]\right\},\cr
&& R_4(\hat{\bbeta}_n^{SL};W)=\sigma^2\tr(WC^{-1})-k\sigma^2\tr(C^{-1}-A)X(\Delta^2)+k(k+4)\bdelta^T W \bdelta E\left[\chi_{q+4}^{-4}(\Delta^2)\right],\cr
&& M_4(\hat{\bbeta}_n^{SL})=\sigma^2 C^{-1}-k\sigma^2 (C^{-1}-A)X(\Delta^2)+k(k+4)\bdelta \bdelta^T E\left[\chi_{q+4}^{-4}(\Delta^2)\right],\cr
\cr
(v) && b_5(\hat{\bbeta}_n^{PRL})=b_4(\hat{\bbeta}_n^{SL})-\bdelta E\left[(1-k\chi_{q+2}^{-2}(\Delta^2))I(\chi_{q+2}^2(\Delta^2)\leq k)\right],\cr
&& B_5(\hat{\bbeta}_n^{PRL})=\Delta^2\left\{kE\left[\chi_{q+2}^{-2}(\Delta^2)\right]-E\left[(1-k\chi_{q+2}^{-2}(\Delta^2))I(\chi_{q+2}^2(\Delta^2)\leq k)\right]\right\}^2,\cr
&& R_5(\hat{\bbeta}_n^{PRL};W)=R_4(\hat{\bbeta}_n^{SL};W)-\sigma^2 \tr(C^{-1}-A)E\left[(1-k\chi_{q+2}^{-2}(\Delta^2))^2I(\chi_{q+2}^2(\Delta^2)\leq k)\right]\cr
&& \hspace*{2.5cm} \quad -\bdelta^T W \bdelta Q(\Delta^2),\cr
&& M_5(\hat{\bbeta}_n^{PRL})=M_4(\hat{\bbeta}_n^{SL})-\sigma^2 (C^{-1}-A)E\left[(1-k\chi_{q+2}^{-2}(\Delta^2))^2I(\chi_{q+2}^2(\Delta^2)\leq k)\right]\cr
&& \hspace*{2.5 cm}\quad -\bdelta\bdelta^T Q(\Delta^2).\nonumber
\end{eqnarray}
where
\begin{eqnarray}
Z(\alpha;\Delta^2)&=& 2H_{q+2}(\chi_q^2(\alpha);\Delta^2)-H_{q+4}(\chi_q^2(\alpha);\Delta^2),\cr
X(\Delta^2)&=& 2E\left[\chi_{q+2}^{-2}(\Delta^2)\right]-k E\left[\chi_{q+4}^{-2}(\Delta^2)\right],\cr
Q(\Delta^2)&=& 2E\left[(1-k\chi_{q+2}^{-2}(\Delta^2))I(\chi_{q+2}^2(\Delta^2)\leq k)\right]-E\left[(1-k\chi_{q+4}^{-2}(\Delta^2))I(\chi_{q+4}^2(\Delta^2)\leq k)\right]\nonumber
\end{eqnarray}
\end{theorem}

\subsection{Graphical Representations}
In this section,  some graphical illustrations will be provided for asymptotic distributional quadratic risk functions. For our purpose, we assume $p=4$, $q=3$,
\begin{equation*}
\bbeta=[1 \ 0 \ -1\ 1]^T,\quad\quad
H=\left[\begin{array}{cccc}
 1 & -1 & 3 & 1 \\
 3 & 2 & 1 &  0\\
 4 & -2 & 0 & 5\\
\end{array}\right],\quad \bh=[0 \ 0 \ 0]^T\quad
 \text{and} \quad \bxi=[1 \ 1 \ 1]^T
\end{equation*}

\begin{figure}
\begin{center}
\begin{tabular}{cc}
\includegraphics[height=7cm,width=7cm]{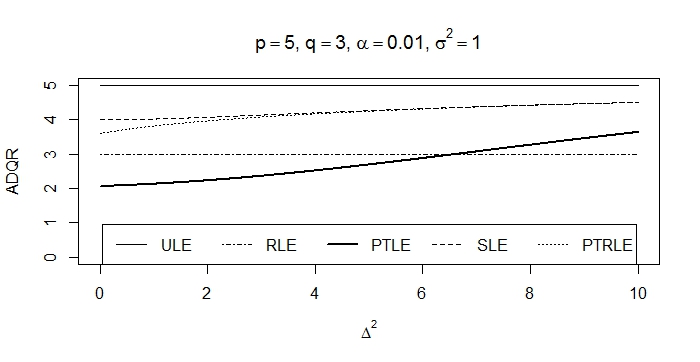}&
\includegraphics[height=7cm,width=7cm]{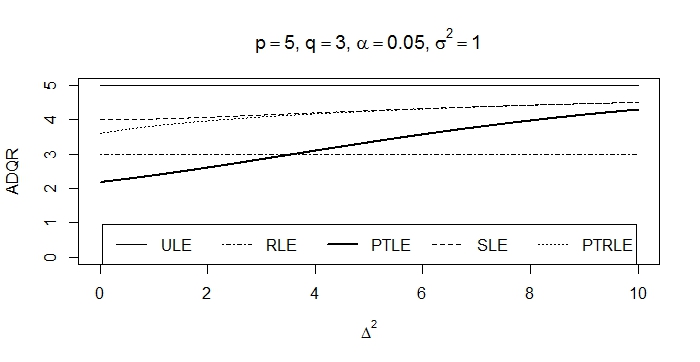}\\
\includegraphics[height=7cm,width=7cm]{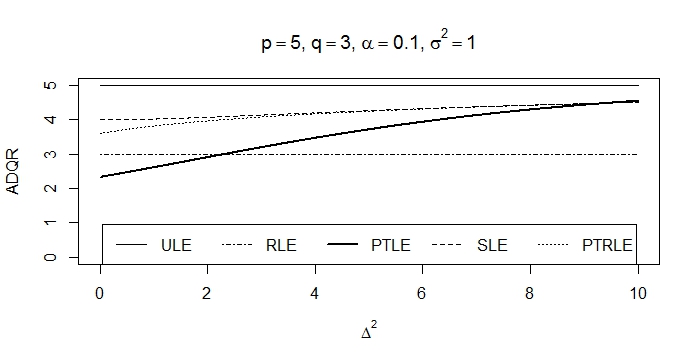}&
\includegraphics[height=7cm,width=7cm]{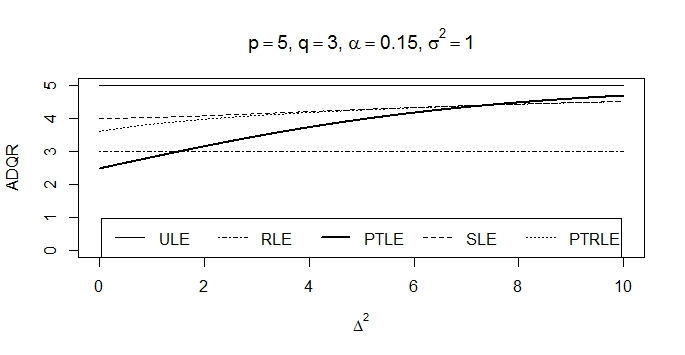}\\
\includegraphics[height=7cm,width=7cm]{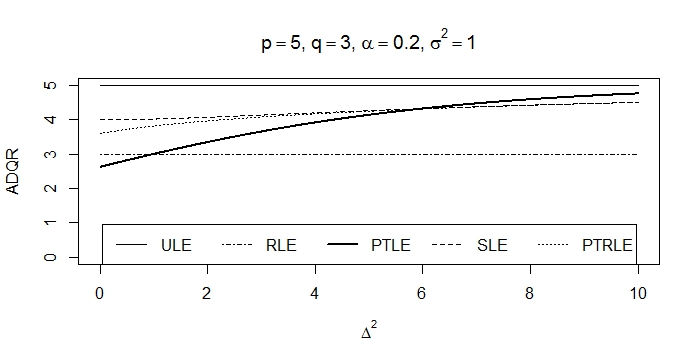}&
\includegraphics[height=7cm,width=7cm]{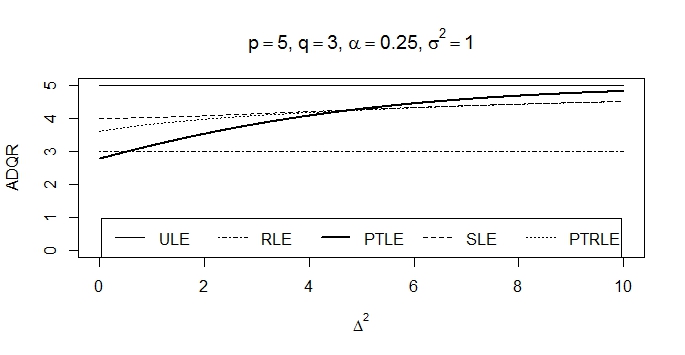}\\
\end{tabular}
\caption{ADQR functions for $\sigma^2=1$, and different level of significance $\alpha$. }
\end{center}
\end{figure}

\begin{figure}
\begin{center}
\begin{tabular}{cc}
\includegraphics[height=7cm,width=7cm]{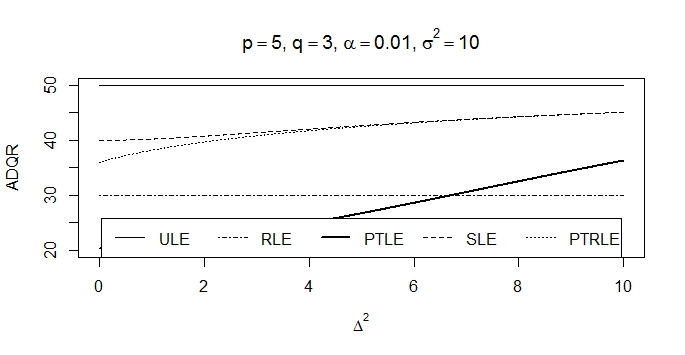}&
\includegraphics[height=7cm,width=7cm]{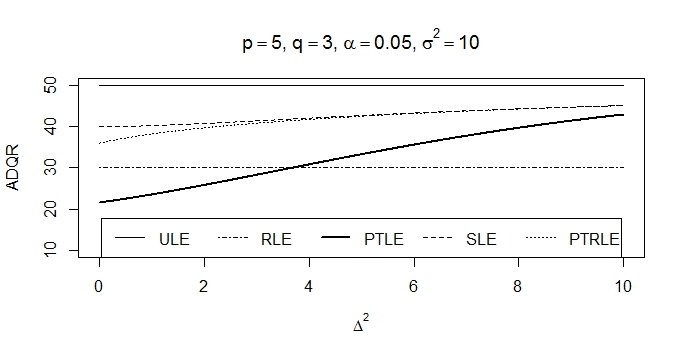}\\
\includegraphics[height=7cm,width=7cm]{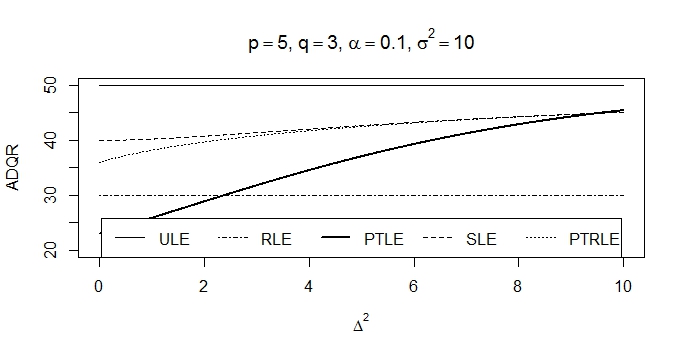}&
\includegraphics[height=7cm,width=7cm]{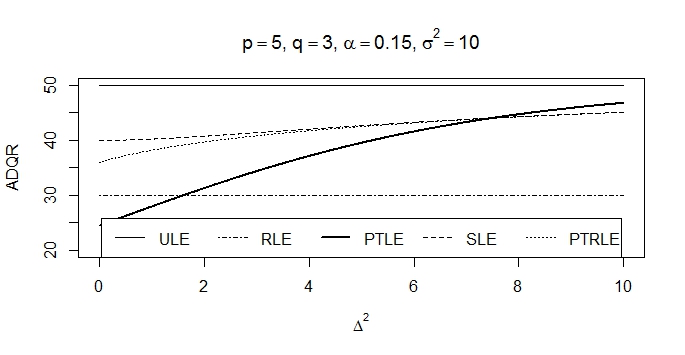}\\
\includegraphics[height=7cm,width=7cm]{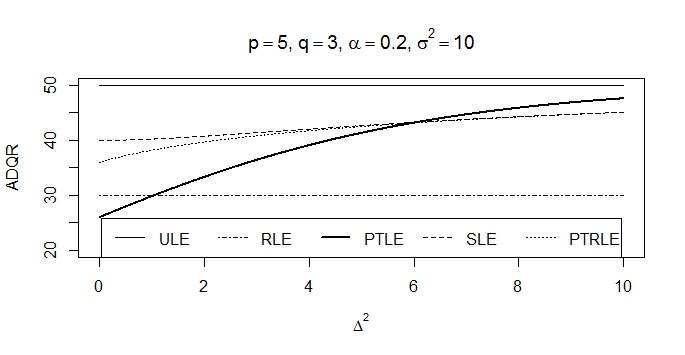}&
\includegraphics[height=7cm,width=7cm]{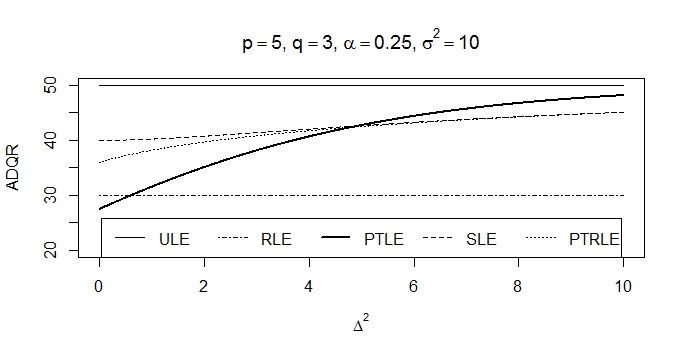}\\
\end{tabular}
\caption{ADQR functions for  $\sigma^2=10$, and different level of significance $\alpha$.}
\end{center}
\end{figure}

From Figures 1 \& 2, it can be deduced that all proposed estimators namely restricted LASSO, preliminary test LASSO, Stein-type Shrinkage LASSO and its positive part estimators perform better than the LASSO estimator in the sense of having smaller ADQR. It is also evident that as we deviate from the null-hypothesis, ADQR values get larger. Finally, as one may expect, by decreasing the level of significant, the preliminary test LASSO estimator performs better.

We also depicted the ADQR functions for different $p$ and $q$ values and found no substantial change in performances.
\newpage
\section{Simulation}\label{sec:4}
In this section, we conduct a Monte Carlo simulation to analyze relative efficiencies with respect to different level of sparsity.

We generate $X$-matrix from a multivariate normal distribution with mean vector $\bmu=0$ and covariance matrix $\Sigma$. The off-diagonal elements of the covariance matrix are considered to be equal to $r$ with $r=0,0.2,0.9$. We consider $n=100$ and various $p$ ranging from $10$ to $30$.

In our simulation scheme, $\bbeta$ is a $p$-vector and a function of $\Delta^2$. When $\Delta^2=0$, $\bbeta$ is the null vector. $\Delta^2>0$ is equivalent to ``violation" of the null hypothesis. We considered $9$ different values for $\Delta^2$, which are $0,1,2,3,5,10,20,30,50$. The way the $\bbeta$ vector is defined in our setup, a $\Delta^2$ indicates that data are generated under null hypothesis, whereas $\Delta^2>0$ indicates a data set generated under alternative hypothesis. Each realization was repeated $2000$ times to obtain bias-squared and variance of the estimated regression parameters.

Finally, risks are calculated for the  ULE, RLE, PTLE, SSLE and PRLE. The responses were simulated from the following model:
\begin{equation*}
y_i=\sum_{i=1}^p X_i\beta_i+e_i, \qquad e_i\sim N(0,5^2)
\end{equation*}
Relative efficiencies are calculated as $\text{Risk}(\hat{\bbeta}^{L})/\text{Risk}(\hat{\bbeta}^*)$, where $\hat{\beta}^*$ is one of the estimators whose relative efficiency is to be computed.

For comparing the relative efficiencies of the penalty estimators, the data generation setup was slightly modified to accomodate the number of nonzero $\beta$s in the model. In particular, we partitioned $\beta$ as $\beta=(k,q)^T$ where $k$ indicates number of nonzero $\beta$s and $q$ indicates $p-k$ zeros-a function of $\Delta^2$. To translate the above, when $p=10$ and $k=5$, we would have $\beta=(1,1,1,1,1,0,0,0,0,0)^T$, and the previously mentioned procedures would be used to generate the data.

From Tables 1-3, it can be realized that the PRSSLE has the best performance among all. As a numerical proof for the assertion in graphical representation, when we deviate from the null model, neither PTLE nor SLE dominates one another and the PTLE performs better as $\alpha$ gets larger. Relative efficiency of the proposed estimators increases when there are more near-zero parameters present in the model. Performance of the estimators decrease as we deviate from the null model.
 \begin{table}[htbp]
 \tiny
 \centering
 \caption{Relative efficiencies of the estimators for fixed $\Delta^2$, $r=0$, different values of $p$ and $k$.}\label{tab:1}
\begin{tabular}{|c|ccccccc|c|ccc|c|c|}
\hline
 & \multicolumn{7}{c|}{ULE} & RLE & \multicolumn{3}{c|}{PTLE} & SSLE & PRLE\\
\hline
 & $ k=1 $ & $ k=3 $ & $ k=4 $ & $ k=5 $ & $ k=6 $ && $ k=p $ &  & $0.15 $ & $0.20$ & $0.25$ & &\\
\hline
$p$ &\multicolumn{13}{c|}{$\Delta^2=0$}\\
\hline
10 & 23.42 & 15.97 & 17.12 & 19.13 & 19.54 & & 19.49 & 96.41 & 24.57 & 22.91 & 21.25 & 39.41 & 39.42\\
20 &43.80 & 46.44 & 47.19 & 43.61 & 47.33& & 44.66 & 347.55 & 64.67 & 58.41 & 50.63 & 144.03 & 144.11\\
30 & 112.44 & 100.43 & 89.37 & 119.88 & 92.40& & 119.22 & 1181.14 & 706.06 & 515.65 & 336.90 & 876.35 & 934.60\\
\hline
$p$ &\multicolumn{13}{c|}{$\Delta^2=1$}\\
\hline
10 & 13.43 & 9.77 & 10.41 & 8.94 & 8.30 & & 6.98 & 9.19 & 7.55 & 7.38 & 7.20 & 8.17 & 8.18\\
20 & 38.07 & 26.55 & 24.00 & 23.79 & 22.79& & 11.09 & 13.42 & 11.72 & 11.56 & 11.33 & 12.88 & 12.88\\
30 & 76.32 & 54.95 & 51.71 & 49.40 & 44.92 & & 17.77 & 20.60 & 20.34 & 20.07 & 19.75 & 20.50 & 20.51\\
\hline
$p$ &\multicolumn{13}{c|}{$\Delta^2=2$}\\
\hline
10 & 10.30 & 8.91 & 7.55 & 7.41 & 6.41 & & 5.53 & 7.13 & 5.97 & 5.81 & 5.68 & 6.41 & 6.41\\
20 & 27.23 & 21.96 & 19.14 & 17.79 & 17.35 & & 10.21 & 12.44 & 10.69 & 10.51 & 10.34 & 11.91 & 11.91\\
30 & 50.46 & 42.27 & 45.57 & 39.07 & 35.93 & & 16.35 & 18.92 & 18.52 & 18.25 & 17.93 & 18.81 & 18.81\\
\hline
$p$ &\multicolumn{13}{c|}{$\Delta^2=3$}\\
\hline
10 & 6.64 & 6.39 & 6.06 & 5.40 & 5.19& & 4.74 & 5.83 & 4.95 & 4.88 & 4.81 & 5.34 & 5.35 \\
20 & 17.34 & 15.74 & 14.55 & 14.00 & 13.44 & & 8.96 & 10.45 & 9.36 & 9.23 & 9.08 & 10.09 & 10.10 \\
30 & 37.58 & 34.27 & 31.82 & 30.83 & 28.71 & & 14.85 & 17.05 & 16.84 & 16.72 & 16.32 & 16.95 & 16.97 \\
\hline
$p$ &\multicolumn{13}{c|}{$\Delta^2=5$}\\
\hline
 10 & 3.97 & 3.77 & 3.50 & 3.58 & 3.38 & & 3.19 & 3.53 & 3.27 & 3.24 & 3.23 & 3.39 & 3.40 \\
 20 & 9.42 & 8.86 & 8.43 & 8.46 & 8.15 & & 6.42 & 7.02 & 6.60 & 6.55 & 6.47 & 6.88 & 6.88 \\
 30 & 20.18 & 18.33 & 17.84 & 18.14 & 17.22 & & 10.97 & 12.02 & 11.99 & 11.89 & 11.70 & 11.99 & 11.99 \\
\hline
$p$ &\multicolumn{13}{c|}{$\Delta^2=10$}\\
\hline
10 & 1.95 & 1.93 & 1.80 & 1.90 & 1.88 & & 1.84 & 1.93 & 1.87 & 1.86 & 1.85 & 1.89 & 1.89 \\
 20 & 3.72 & 3.66 & 3.68 & 3.65 & 3.55  & & 3.33 & 3.54 & 3.38 & 3.36 & 3.34 & 3.49 & 3.49 \\
 30 & 6.93 & 7.09 & 6.99 & 6.87 & 6.88 & & 6.18 & 6.53 & 6.49 & 6.47 & 6.43 & 6.51 & 6.51 \\
 \hline
 $p$ &\multicolumn{13}{c|}{$\Delta^2=20$}\\
 \hline
10 & 1.35 & 1.32 & 1.35 & 1.33 & 1.36 & & 1.36 & 1.38 & 1.37 & 1.36 & 1.36 & 1.37 & 1.37 \\
 20 & 1.99 & 2.01 & 2.01 & 1.97 & 2.00 & & 1.97 & 2.00 & 1.98 & 1.98 & 1.97 & 1.99 & 1.99 \\
 30 &  3.22 & 3.25 & 3.21 & 3.25 & 3.30  & & 3.25 & 3.32 & 3.31 & 3.31 & 3.30 & 3.32 & 3.32 \\
 \hline
  $p$ &\multicolumn{13}{c|}{$\Delta^2=30$}\\
  \hline
 10 &1.23 & 1.22 & 1.21 & 1.23 & 1.22 & & 1.23 & 1.24 & 1.24 & 1.23 & 1.23 & 1.24 & 1.24 \\
  20 & 1.63 & 1.65 & 1.64 & 1.66 & 1.61 & & 1.68 & 1.71 & 1.69 & 1.69 & 1.69 & 1.70 & 1.70 \\
  30 & 2.47 & 2.54 & 2.54 & 2.52 & 2.52 & & 2.48 & 2.51 & 2.51 & 2.50 & 2.50 & 2.51 & 2.51 \\
 \hline
  $p$ &\multicolumn{13}{c|}{$\Delta^2=50$}\\
  \hline
 10 & 1.17 & 1.17 & 1.18 & 1.16 & 1.18  & &  1.16 & 1.17 & 1.17 & 1.17 & 1.17 & 1.17 & 1.17 \\
  20 & 1.49 & 1.53 & 1.47 & 1.50 & 1.49 & & 1.51 & 1.53 & 1.52 & 1.52 & 1.51 & 1.52 & 1.52 \\
  30 & 2.14 & 2.13 & 2.08 & 2.11 & 2.09 & & 2.15 & 2.18 & 2.17 & 2.17 & 2.17 & 2.18 & 2.18 \\
  \hline
\end{tabular}
 \end{table}

  \begin{table}[htbp]
  \tiny
  \centering
  \caption{Relative efficiencies of the estimators for fixed $\Delta^2$, $r=0.2$, different values of $p$ and $k$.}\label{tab:2}
 \begin{tabular}{|c|ccccccc|c|ccc|c|c|}
 \hline
  & \multicolumn{7}{c|}{ULE} & RLE & \multicolumn{3}{c|}{PTLE} & SSLE & PRLE\\
 \hline
  & $ k=1 $ &$  k=3 $ & $ k=4 $ &$  k=5 $ & $ k=6  $&&$  k=p $ &  & $0.15 $ & $0.20$ & $0.25$ & &\\
 \hline
 $p$ &\multicolumn{13}{c|}{$\Delta^2=0$}\\
 \hline
 10 & 20.64 & 20.21 & 16.21 & 18.22 & 20.23 & & 18.14 & 73.26 & 23.72 & 21.20 & 19.70 & 34.63 & 34.91 \\
 20 & 44.37 & 47.21 & 45.70 & 53.74 & 44.38 & & 50.79 & 389.22 & 65.82 & 59.01 & 54.34 & 162.41 & 163.14 \\
 30 & 108.61 & 120.92 & 78.55 & 116.06 & 86.10 & & 114.17 & 673.36 & 477.85 & 442.09 & 289.67 & 576.68 & 620.72 \\
  \hline
 $p$ &\multicolumn{13}{c|}{$\Delta^2=1$}\\
 \hline
 10 & 14.67 & 12.02 & 9.59 & 9.72 & 9.03 & & 7.38 & 10.33 & 8.02 & 7.77 & 7.65 & 8.98 & 9.02 \\
 20 & 35.14 & 29.74 & 26.11 & 27.71 & 23.27& & 13.62 & 17.33 & 14.75 & 14.30 & 14.05 & 16.54 & 16.54 \\
 30 & 83.94 & 72.43 & 51.38 & 59.23 & 46.21& & 24.32 & 29.79 & 29.17 & 28.94 & 27.99 & 29.55 & 29.68 \\
 \hline
 $p$ &\multicolumn{13}{c|}{$\Delta^2=2$}\\
 \hline
 10 &  11.13 & 8.78 & 8.72 & 8.53 & 7.66 & & 6.29 & 8.48 & 6.77 & 6.58 & 6.46 & 7.45 & 7.45\\
 20 & 30.16 & 25.60 & 22.51 & 21.76 & 21.15& &12.71 & 15.62 & 13.68 & 13.30 & 12.90 & 14.99 & 15.00 \\
 30 & 63.08 & 49.59 & 51.37 & 45.83 & 43.62 & & 23.52 & 28.94 & 27.87 & 27.32 & 26.72 & 28.62 & 28.70 \\
 \hline
 $p$ &\multicolumn{13}{c|}{$\Delta^2=3$}\\
 \hline
 10 &  7.19 & 6.77 & 6.55 & 6.38 & 6.25 & & 5.75 & 6.69 & 6.03 & 5.93 & 5.85 & 6.30 & 6.31 \\
 20 & 19.67 & 17.27 & 17.39 & 15.67 & 15.72 & &  11.02 & 13.38 & 11.56 & 11.43 & 11.19 & 12.86 & 12.86 \\
 30 & 42.76 & 40.57 & 40.57 & 34.44 & 33.20 & &  20.80 & 24.50 & 24.13 & 23.93 & 23.04 & 24.33 & 24.40 \\
 \hline
 $p$ &\multicolumn{13}{c|}{$\Delta^2=5$}\\
 \hline
  10 & 4.38 & 4.15 & 4.09 & 3.80 & 3.96 && 3.65 & 4.02 & 3.75 & 3.72 & 3.69 & 3.88 & 3.88 \\
  20 & 10.96 & 10.39 & 9.77 & 10.02 & 9.70 && 8.36 & 9.45 & 8.67 & 8.57 & 8.48 & 9.23 & 9.23 \\
  30 & 23.94 & 21.94 & 20.90 & 20.14 & 20.15 && 17.11 & 19.63 & 19.29 & 19.17 & 18.66 & 19.48 & 19.54 \\
 \hline
 $p$ &\multicolumn{13}{c|}{$\Delta^2=10$}\\
 \hline
 10 & 2.03 & 2.02 & 2.03 & 2.09 & 2.06 && 2.03 & 2.11 & 2.06 & 2.05 & 2.04 & 2.08 & 2.08 \\
  20 & 4.12 & 4.08 & 4.22 & 4.05 & 4.28 && 4.34 & 4.61 & 4.42 & 4.40 & 4.37 & 4.56 & 4.56 \\
  30 &   8.63 & 8.35 & 8.40 & 8.49 & 8.43 && 9.28 & 9.91 & 9.78 & 9.74 & 9.62 & 9.87 & 9.88 \\
  \hline
  $p$ &\multicolumn{13}{c|}{$\Delta^2=20$}\\
  \hline
 10 &  1.41 & 1.38 & 1.43 & 1.41 & 1.44  && 1.42 & 1.44 & 1.42 & 1.42 & 1.42 & 1.43 & 1.43 \\
  20 & 2.15 & 2.21 & 2.24 & 2.19 & 2.25 && 2.45 & 2.50 & 2.46 & 2.46 & 2.46 & 2.49 & 2.49 \\
  30 &  3.85 & 3.74 & 3.64 & 3.91 & 3.79 && 4.61 & 4.76 & 4.74 & 4.74 & 4.71 & 4.75 & 4.75 \\
  \hline
   $p$ &\multicolumn{13}{c|}{$\Delta^2=30$}\\
   \hline
  10 & 1.28 & 1.26 & 1.25 & 1.26 & 1.28 && 1.42 & 1.44 & 1.42 & 1.42 & 1.42 & 1.43 & 1.43 \\
   20 & 1.83 & 1.82 & 1.77 & 1.79 & 1.79 && 2.45 & 2.50 & 2.46 & 2.46 & 2.46 & 2.49 & 2.49 \\
   30 & 2.83 & 2.82 & 2.91 & 2.89 & 2.90 && 4.61 & 4.76 & 4.74 & 4.74 & 4.71 & 4.75 & 4.75 \\
  \hline
   $p$ &\multicolumn{13}{c|}{$\Delta^2=50$}\\
   \hline
  10 &  1.19 & 1.18 & 1.20 & 1.21 & 1.20 && 1.21 & 1.22 & 1.22 & 1.21 & 1.21 & 1.22 & 1.22 \\
   20 & 1.60 & 1.58 & 1.63 & 1.61 & 1.64 && 1.65 & 1.67 & 1.66 & 1.65 & 1.65 & 1.67 & 1.67 \\
   30 &  2.36 & 2.42 & 2.47 & 2.45 & 2.41 && 2.65 & 2.69 & 2.68 & 2.68 & 2.67 & 2.69 & 2.69 \\
   \hline
 \end{tabular}
  \end{table}

    \begin{table}[htbp]
    \tiny
    \centering
    \caption{Relative efficiencies of the estimators for fixed $\Delta^2$, $r=0.9$, different values of $p$ and $k$.}\label{tab:3}
   \begin{tabular}{|c|ccccccc|c|ccc|c|c|}
   \hline
    & \multicolumn{7}{c|}{ULE} & RLE & \multicolumn{3}{c|}{PTLE} & SSLE & PRLE\\
   \hline
    & $ k=1 $ &$  k=3 $ &$  k=4  $& $ k=5 $ &$k=6$ && $k=p$ &  & $0.15 $ & $0.20$ & $0.25$ & &\\
   \hline
   $p$ &\multicolumn{13}{c|}{$\Delta^2=0$}\\
   \hline
    10 & 20.82 & 20.22 & 16.11 & 19.35 & 21.15 && 19.03 & 123.65 & 23.36 & 20.80 & 20.00 & 38.20 & 38.74 \\
    20  & 45.05 & 46.89 & 46.75 & 53.90 & 44.56 && 50.71 & 632.12 & 67.21 & 60.49 & 53.35 & 173.14 & 175.40 \\
    30 & 108.96 & 118.30 & 79.47 & 115.30 & 85.65 && 114.00 & 1172.26 & 759.67 & 585.21 & 317.76 & 826.66 & 912.68 \\
    \hline
   $p$ &\multicolumn{13}{c|}{$\Delta^2=1$}\\
   \hline
  10 & 21.52 & 17.79 & 21.27 & 21.48 & 17.66 & & 15.04 & 53.45 & 18.09 & 16.77 & 15.79 & 25.53 & 26.14 \\
  20 & 53.00 & 59.20 & 47.68 & 49.57 & 47.64 & & 42.22 & 140.71 & 51.69 & 45.55 & 43.08 & 94.46 & 94.96 \\
   30 & 104.08 & 92.84 & 98.95 & 96.73 & 90.50 & & 85.50 & 344.20 & 226.35 & 208.81 & 172.24 & 298.30 & 306.38 \\
   \hline
   $p$ &\multicolumn{13}{c|}{$\Delta^2=2$}\\
   \hline
   10 & 15.59 & 15.52 & 16.25 & 16.19 & 19.80 &   & 21.24 & 52.47 & 26.07 & 22.69 & 21.51 & 32.62 & 33.01 \\
   20 & 47.47 & 42.76 & 46.90 & 39.31 & 47.51 &   & 40.03 & 134.72 & 47.79 & 44.32 & 41.59 & 89.99 & 90.42 \\
   30 & 97.37 & 105.75 & 106.75 & 88.70 & 82.06 &   & 95.76 & 308.65 & 264.45 & 222.31 & 165.95 & 275.59 & 285.19 \\
   \hline
   $p$ &\multicolumn{13}{c|}{$\Delta^2=3$}\\
   \hline
   10& 19.50 & 16.80 & 15.94 & 14.38 & 15.77   & & 16.35 & 36.63 & 18.22 & 17.82 & 17.25 & 24.47 & 24.82 \\
    20& 42.78 & 44.65 & 40.20 & 44.48 & 37.82   & & 41.41 & 118.89 & 50.64 & 45.98 & 43.02 & 86.71 & 87.32 \\
    30& 99.68 & 90.81 & 82.79 & 73.58 & 74.62   & & 83.46 & 297.35 & 218.31 & 198.30 & 143.54 & 260.92 & 268.40 \\
   \hline
   $p$ &\multicolumn{13}{c|}{$\Delta^2=5$}\\
   \hline
  10& 12.24 & 12.50 & 11.49 & 12.89 & 12.51 & & 12.38 & 24.72 & 14.18 & 13.39 & 12.86 & 17.75 & 17.81 \\
  20& 34.70 & 28.45 & 32.99 & 34.51 & 35.30  & & 32.62 & 86.57 & 41.94 & 38.01 & 35.64 & 65.12 & 65.35 \\
   30&71.34 & 61.56 & 73.58 & 64.65 & 78.40  & & 76.18 & 234.67 & 190.59 & 156.13 & 123.91 & 212.46 & 217.32 \\
   \hline
   $p$ &\multicolumn{13}{c|}{$\Delta^2=10$}\\
   \hline
  10& 7.72 & 7.14 & 7.35 & 7.88 & 7.83 & & 8.10 & 11.70 & 8.85 & 8.66 & 8.33 & 9.96 & 9.99 \\
   20& 19.24 & 19.80 & 18.82 & 20.82 & 20.12 & & 24.23 & 42.53 & 26.37 & 25.72 & 24.70 & 36.54 & 36.62 \\
   30& 41.44 & 43.91 & 41.64 & 44.13 & 43.36 & & 54.70 & 125.52 & 102.55 & 89.92 & 80.92 & 116.12 & 117.89 \\
    \hline
    $p$ &\multicolumn{13}{c|}{$\Delta^2=20$}\\
    \hline
  10& 3.70 & 3.93 & 4.00 & 4.08 & 4.30   & & 4.33 & 5.19 & 4.50 & 4.44 & 4.34 & 4.81 & 4.82 \\
  20& 9.77 & 9.67 & 9.90 & 10.32 & 10.62 & & 14.67 & 19.13 & 15.67 & 15.46 & 15.07 & 17.95 & 17.98 \\
  30& 19.07 & 20.16 & 21.57 & 22.40 & 22.48 & &37.61 & 56.97 & 54.06 & 51.10 & 47.37 & 55.21 & 55.64\\
    \hline
     $p$ &\multicolumn{13}{c|}{$\Delta^2=30$}\\
     \hline
    10 & 2.99 & 2.99 & 3.19 & 3.28 & 3.30 && 3.48 & 3.93 & 3.55 & 3.54 & 3.49 & 3.70 & 3.71 \\
     20 & 6.85 & 7.27 & 7.43 & 7.55 & 7.90 && 10.44 & 12.56 & 11.00 & 10.85 & 10.61 & 11.99 & 12.03 \\
     30 & 14.29 & 15.10 & 15.41 & 16.12 & 16.41 & & 27.59 & 36.70 & 35.27 & 34.41 & 32.41 & 36.34 & 36.37 \\
    \hline
     $p$ &\multicolumn{13}{c|}{$\Delta^2=50$}\\
     \hline
  10 & 2.46 & 2.51 & 2.52 & 2.61 & 2.69 & & 2.70 & 2.94 & 2.75 & 2.73 & 2.71 & 2.82 & 2.83 \\
  20&  5.53 & 5.64 & 5.62 & 5.92 & 6.08 & & 7.68 & 8.74 & 7.95 & 7.84 & 7.77 & 8.50 & 8.51 \\
  30& 11.07 & 11.86 & 11.72 & 12.27 & 12.80 & & 18.75 & 23.69 & 22.71 & 22.14 & 21.81 & 23.42 & 23.44 \\
     \hline
   \end{tabular}
    \end{table}

\newpage
\section{Real Data}\label{sec:5}
\begin{figure}
\centering
\includegraphics[scale=0.5]{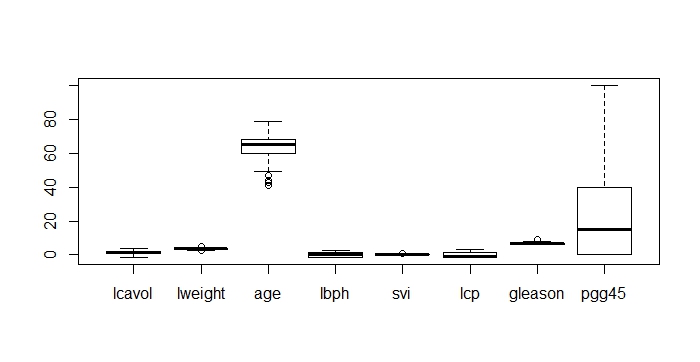}
\caption{The box plot of predictors in prostate data.}\label{fig:5}
\end{figure}

In this section we study the performance of the proposed LASSO-type estimators in a real example. We use the prostate dataset (Stamey et al, 1989). These data come from a study that examined the correlation between the level of prostate specific antigen and a number of clinical measures in men who were about to receive a radical prostatectomy. A descriptions of  the variables in this dataset is given in Table \ref{tab:4} and Figure \ref{fig:5} shows the box-plot of the variables.

\begin{table}[htnp]
\small
\centering
\caption{Discription of the variables of prostate data. }\label{tab:4}
\begin{tabular}{lll}
\hline
Variables & Description &  Remarks\\
\hline
lpsa & Log of prostate specific antigen (PSA) & Response\\
lcavol & Log cancer volume & \\
lweight & Log prostate weight &\\
age & Age  & \\  Age in years.
lbph & Log of benign prostatic hyperplasia amount & \\
svi & Seminal vesicle invasion & \\
lcp & Log of capsular penetration & \\
gleason &  Gleason score& A  numeric vector\\
pgg45 & Percent of Gleason scores 4 or 5& \\
\hline
\end{tabular}
\end{table}

First, we center the predictor variables. At the second step, we fit the linear regression model to predict the response variable in the presence of regressors. Unrestricted LASSO (UL), restricted LASSO (RL), preliminary test LASSO (PTL), Stein-type shrinkage LASSO (SL), and positive rule Stein-type shrinkage (PRL) estimators are used to estimate the regression parameters.

The summary statistics of response variable (lpsa) is shown in Table \ref{tab:5}.
\begin{table}
\centering
\caption{Summary statistics for response variable in the prostate dataset.}\label{tab:5}
\begin{tabular}{ccccccc}
\hline
Min & Q1 & Median & Q3 & Max &  Mean  & SD \\
\hline
$-0.4308$ &$1.7320$ & $2.5920$ & $3.0560$ & $5.5830$ & $2.478$0 & $1.1543$\\
\hline
\end{tabular}
\end{table}

The performance of the estimators are evaluated using average $10$-fold cross validation error. $k$-fold cross validation is a famous method that divide the data set into $k$ equal-seized subset, randomly. One of the subsets is selected as test set and the $k-1$ renaming subsets are called train set and used to fit the model. The obtained model is then used for predicting the response variable in the test set. Prediction errors is the squared version of difference between the observed and predicted values of the response variable in  the test set.

We used the following specifications:
\begin{equation*}
H=\begin{bmatrix}
$-1$&$3$&$1$&$-1$&$0$&$-1$&$0$&$0$\\
$-1$&$1$&$0$&$-1$&$0$&$1$&$0$&$0$\\
$1$&$0$&$-1$&$1$&$0$&$0$&$1$&$0$
\end{bmatrix}, \qquad h=\begin{bmatrix}
0\\ 0\\ 0
\end{bmatrix}.
\end{equation*}

By choosing $1000$ as a large enough number for repeating process in a bootstrap simulation scheme, Table \ref{tab:6} shows the average and standard deviation of the prediction errors.

\begin{table}[ht]
\centering
\caption{$10$-fold cross validation average prediction errors and standard deviations for prostate data}\label{tab:6}
\begin{tabular}{lccccccc}
  \hline
 & ULE & RLE & PTLE(0.01) & PTLE(0.05) & PTLE(0.10) & SSLE & PRSSLE \\
  \hline
mean & 5822.99 & 5753.18 & 6122.17 & 6076.60 & 5906.45 & 5798.20 & 5797.90 \\
  sd & 92.56 & 69.58 & 251.41 & 262.18 & 213.98 & 92.43 & 92.54 \\
   \hline
\end{tabular}
\end{table}

Based on Table \ref{tab:6}, PRSSLE is the best estimator in terms of the prediction error (the lesser risk, the better estimator). This estimator is followed by SSLE.  If the level of significance $\alpha$ for constructing PTLE increases, then the prediction error decreases. These results confirm our assertions.

\begin{figure}
\centering
\includegraphics[scale=0.7]{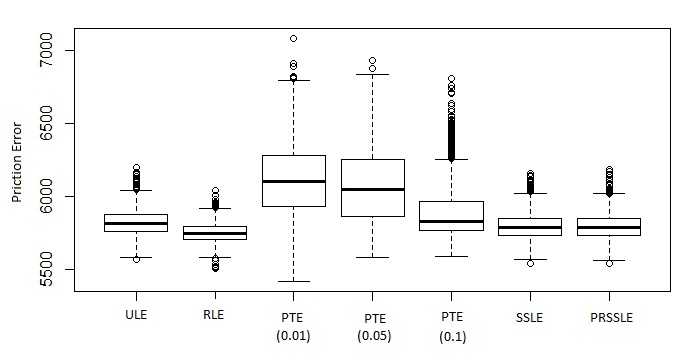}
\caption{Boxplot of prediction errors of the unrestricted LASSO (UL), restricted LASSO (RL), preliminary-test LASSO (PTL) for $\alpha=0.01,0.05$ and $0.1$,Stein-type shrinkage LASSO (SSL) and positive rule Stein-type shrinkage LASSO (PRSSL) estimators for prostate data.}\label{fig:6}
\end{figure}

Figure \ref{fig:6} shows the boxplot of average prediction errors for the proposed estimators visually, which demonstrates small prediction errors for the PRSSLE and SSLE. Although PRSSLE has the smallest prediction error, SSLE has the less variability of the five estimators. Variability of RLE is less than others, which means the null hypothesis is approximately true, here.

\section{Conclusion}
In this paper, we proposed an improvement on the LASSO estimator by imposing a restriction to the model and using preliminary testing and shrinkage techniques. Indeed, we introduced preliminary-test LASSO (PRL), Stein-type shrinkage LASSO (SSL) and positive-rule shrinkage LASSO (PRSSL) estimators in the presence of a sub-space restriction.
Performance of the proposed estimators under the null hypothesis has been studied in case of sample size ($n$) is more than the number of features ($p$). The proposed methodology for improving the LASSO can be applied in high dimensional situations, i.e., large $p$, small $n$.

In addition to the given theorems for asymptotic behaviour of the proposed estimators analytically, using a  simulation study, we compared the performance of estimators numerically for various configurations of parameter ($p$), correlation coefficient between the predictors ($r$), and the error in variance ($\sigma^2$). In order to compare the improved estimators with the classical LASSO, we used the relative efficiency criterion.
For different non-centrality parameter $\Delta^2$, degree of model misspecification, the number of non-zero $\bbeta$s varied, and then the performance of estimators evaluated.
We found out that the PRSSLE has the best performance among all. When we deviate from the null model, neither PTLE nor SLE dominates one another and the PTLE performs better as $\alpha$ gets larger. Relative efficiency of the proposed estimators increases when there are more near-zero parameters present in the model. Performance of the estimators decrease as we deviate from the null model.

As an application, the prostate dataset analyzed. In this respect, $10$-folded cross validation average and standard deviations of the prediction errors based on the LASSO , Restricted LASSO, preliminary test LASSO, Stein-type shrinkage LASSO and positive rule Stein type shrinkage LASSO compared. The new estimators dominate the LASSO one in average prediction error sense but the picture of dominance of the PRSSLE was not obvious.  We conclude that the improved LASSO estimators are better than the classical version.
\appendix

\section{Proof of Theorems}

\subsection{Proof of Theorem \ref{thm:4}}
For (i), see Knight and Fu (2000). To prove (ii), by \textit{Sluskey's theorem},  equation \eqref{eq:RLASSO} and also assumption {\bf A2}, we know
\begin{equation*}
\hat{\bbeta}_n^{L}-C_n^{-1}H^T(HC_n^{-1}H^T)^{-1}(H\hat{\bbeta}^{L}_n-\bh)\overset{\mathcal{P}}{\to} \text{argmin}(Z)-C^{-1}H^T(H C^{-1}H^T)^{-1}(H\text{argmin}(Z)-\bh)
\end{equation*}
(iii) By Eq. \eqref{eq:RLASSO}, we have  $\hat{\bbeta}_n^{L}-\hat{\bbeta}_n^{RL}=C_n H^T(H C_n^{-1} H^T)^{-1}(H \hat{\bbeta}_n^{L}-\bh)$, which converges to $C^{-1}H^T (HC^{-1}H^T)^{-1}(H \text{argmin}(Z)-\bh)$ again, by Sluskey's theorem and also assumption {\bf A2}\\
(iv) Base on Theorem \ref{thm:3}, we know $I(\L_{n}\leq \L_{n,\alpha})\overset{\mathcal{D}}{\to} I(\L \leq \L_{\alpha})$. Now, by Eq. \eqref{eq:PTLE}, part (iii) of this theorem and again by Sluskey's  theorem, we have
\begin{equation*}
\hat{\bbeta}_n^{PTL}\overset{\mathcal{D}}{\to}\text{argmin}(Z)-C^{-1}H^T (HC^{-1}H^T)^{-1}(H \text{argmin}(Z)-\bh)I(\L \leq \L_{\alpha})
\end{equation*}
 To prove (v) and (vi),  consider that $k_n\overset{\mathcal{D}}{\to} k=q-2$, based on equation \eqref{eq:SLE}, part (iii) of this theorem and by Sluskey's theorem, the result is obvious.

  \subsection{Proof for Theorem \ref{thm:6}}
  Referring to Sen and Singer (1993) and Sluskey's theorem, (i) is obvious.
  The proof of (ii) is in Knight and Fu (2000). For (iii), by equation \eqref{eq:RLASSO},
  \begin{eqnarray}
  \sqrt{n}(\hat{\bbeta}^{RL}_n-\bbeta)&=&\sqrt{n}(\hat{\bbeta}^{L}_n-\bbeta)-\sqrt{n}C_n^{-1}H^T(H C_n^{-1} H^T)^{-1}[H(\hat{\bbeta}^{L}_n-\bbeta)+(H\bbeta-\bh)]\cr
  &=&\sqrt{n}(\hat{\bbeta}^{L}_n-\bbeta)-C_n^{-1}H^T(H C_n^{-1} H^T)^{-1}[H\sqrt{n}(\hat{\bbeta}^{L}_n-\bbeta)+\sqrt{n}(H\bbeta-\bh)]\cr
  &=&\sqrt{n}(\hat{\bbeta}^{L}_n-\bbeta)-nC_n^{-1}H^T(nH C_n^{-1} H^T)^{-1}[H\sqrt{n}(\hat{\bbeta}^{L}_n-\bbeta)+\sqrt{n}(H\bbeta-\bh)],\nonumber
  \end{eqnarray}
  Making use of Sluskey's theorem and part (ii), we obtain
  \begin{equation*}
  \sqrt{n}(\hat{\bbeta}^{RL}_n-\bbeta)\to_d \text{argmin}(V)-C^{-1}H^T(H C^{-1} H^T)^{-1}[H\text{argmin}(V)+\bxi].
  \end{equation*}
  To prove (iv), we have
  \begin{eqnarray}
  \sqrt{n}(\hat{\bbeta}_n^{L}-\hat{\bbeta}^{RL}_n)&=&\sqrt{n}(\hat{\bbeta}_n^{L}-\hat{\bbeta}^{LE}_n+C_n^{-1}H^T(H C_n^{-1}H^T)^{-1}(H\hat{\bbeta}^{L}_n-\bh))\cr
  &=& \sqrt{n}C_n^{-1}H^T(H C_n^{-1}H^T)^{-1}(H\hat{\bbeta}^{L}_n-\bh)\cr
  &=& C_n^{-1}H^T(H C_n^{-1}H^T)^{-1}(H\sqrt{n}(\hat{\bbeta}^{L}_n-\bbeta)+\sqrt{n}(H\bbeta-\bh))\cr
  &=& C_n^{-1}H^T(H C_n^{-1}H^T)^{-1}(H\sqrt{n}(\hat{\bbeta}^{L}_n-\bbeta)+\bxi).\nonumber
  \end{eqnarray}
  The result is followed by Sluskey's theorem and part (ii).\\
  Part (v) is an immediate consequence of Theorem \ref{thm:3}.
  To prove (vi), by Theorem \ref{thm:3}, using the fact that
  $I(\L_{n}\leq \L_{n,\alpha})\overset{\mathcal{D}}{\to} I(\L \leq \L_{\alpha})$,
   parts (ii) and (iv) and Sluskey's theorem, $\sqrt{n}(\hat{\bbeta}^{PTL}-\bbeta)$ converges in distribution to
  \begin{equation*}
  \text{argmin(V)}-C^{-1}H^T(H C^{-1}H^T)^{-1}(H \text{argmin} V+\bxi)I(\L \leq \L_{\alpha})
  \end{equation*}
  To prove (vii), we have that
  \begin{equation*}
  \sqrt{n}(\hat{\bbeta}^{SSL}-\bbeta)=\sqrt{n}(\hat{\bbeta}^{L}_n-\bbeta)-k_n\sqrt{n}(\hat{\bbeta}^{L}_n-\hat{\bbeta}^{RL}_n)\L_n^{-1}\label{eq:SL-proof}
  \end{equation*}
  By parts (ii) and (iv), $\L_n^{-1}\overset{\mathcal{D}}{\to} \L$, $k_n\overset{\mathcal{D}}{\to} k$ and applying Sluskey's theorem, we may write
  \begin{equation}
  \sqrt{n}(\hat{\bbeta}^{SSL}-\bbeta)\to_d \text{argmin(V)}-k \left[(C^{-1}H^T(H C^{-1}H^T)^{-1}(H \text{argmin} (V)+\bxi)\right]\L^{-1}
  \end{equation}
  And finally for proving (viii), we have
  \begin{equation*}
  \sqrt{n}(\hat{\bbeta}^{PRSL}-\bbeta)=\sqrt{n}(\hat{\bbeta}^{RL}_n-\bbeta)+\left\{k_n\L^{-1}+(1-k_n\L_n^{-1})I(\L_n<k_n)\right\}\sqrt{n}(\hat{\bbeta}^{L}_n-\hat{\bbeta}^{RL}_n)
  \end{equation*}
  In the same fashion as in (vii), and $I(\L_{n}\leq \L_{n,\alpha})\overset{\mathcal{D}}{\to} I(\L \leq \L_{\alpha})$, the proof for $\sqrt{n}(\hat{\bbeta}^{PRSL}-\bbeta)$ is straight.

\subsection{Proof of Theorem \ref{thm:5}}
If $\lambda_n=o(n)$, then $\frac{\lambda_n}{n}\to 0$, i.e. $\lambda_0=0$. We have the minimum of $Z(\phi)$, defined in part (i) of Theorem \ref{thm:4}, as
\begin{equation*}
\frac{\partial}{\partial\phi}Z(\phi)=2c\phi-2c\bbeta=0
\end{equation*}
 It concludes that $\hat{\bbeta}_n^L\overset{\mathcal{P}}{\to}\bbeta$. Thus, $\text{argmin}Z(\phi)=\bbeta$. Implementing this result in Theorem \ref{thm:4}, the proof is complete, since under the null hypothesis $\mathcal{H}_o$, $H\bbeta-\bh=0$ and all estimators are consistent.

\subsection{Proof of Theorem \ref{thm:7}}
In a similar fashion as in the proof of Theorem 7.8.2.3 of Saleh (2006), and using the fact that under $\sqrt{n}$-consistency, $\lambda=0$, $\text{argmin}(V)=W$, all the given results are followed directly, after some algebra.

\section*{References}
\baselineskip=12pt
\def\ref{\noindent\hangindent 25pt}

\ref  Bancroft, T. A. (1944). On biases in estimation due to the use of preliminary tests of significance, {\em Ann. Math. Statist.}, {\bf 15}, 190-204.

\ref Breiman, L. (1996).  Heuristics of instability and stabilization in model selection,  {\em Ann. Statist.}, {\bf 24}(6), 2350--2383.

%\ref Efron, B. (2004). The estimation of prediction error: covariance penalties and cross-validation, {\em J. Amer. Statist. Assoc.}, {\bf 99}(467), 619–642.

\ref Fan, J. and Li, R. (2001). Variable selection via nonconcave penalized likelihood and its oracle properties, {\em J. Amer. Statist. Assoc.}, {\bf 96}(456), 1348–-1360.

\ref Fan, J. Zhang, C. Zhang, J. (2001). Generalized likelihood ratio statistics and wilks phenomenon, {\em Ann. Statist.},{\bf 29}, 153--193.

\ref Gibbons. D. G (1981). A simulation study of some ridge estimators. {\em J. Amer. Statist.Assc.},{\bf 76}, 1131--19.

\ref Gut, A. (2005). {\em Probability: a graduate study}, Springer, The United Stated.

\ref Hastie, T., Tibshirani, R. and Friedman, J. (2009), {\em The Elements of Statistical Learning; Data Mining, Inference and Prediction}, Springer Verlag, New York.

\ref Hoerl. A. E., Kennard. R. W.(1970). Ridge regression baised estimation for non-orthogonal problems. {\em Thechnometrics}, {\bf 12},69--89.

\ref Knight, K. and Fu, W. (2000). Asymptotics for LASSO-type estimators, {\em Annals of Statistics}, {\bf 28},(5), 1356--1378.
\ref  Mc.Donald.G.C, Galarneau.D.I. (1975). A monte carlo evalution of some ridge type estimators.{\em J. Amer. Statist. Assoc.},{\bf 70},407--416.

\ref Miller, A. J. (2002). {\em Subset selection in Regression}, 2$^{nd}$ Ed., Chapman \& Hall\/CRC, USA.

%\ref Rao.C.R.(1995).{\em Linear Models-Least Squares and Alternatives}, Springer.
%\ref Kunugi, T., Tamura, T., Natio, T. (1961) New acetylene process uses hydrogen dilution.{\em Chen. Eng. Prog. }, {\bf 57}, 43--49.

\ref Saleh, A. K. Md. Ehsanes, (2006). Theory of Preliminary Test and Stein-Type Estimation with Applications, Wiley; United Stated of America.

\ref Sen, P. K. and Singer, J. M. (1993). {\em Large Sample Methods in Statistics: An Introduction with applications}, Chapman and Hall, New York.

\ref Sengupta, D. and Jammalamadaka, S. R. (2003). {\em Linear Models: An Integrated Approach}, World Scientific Publishing Company, Singapore.

\ref Shorack, ?. (2000).{\em Probability for Statisticians}, Springer Texts in Statistics. Springer-Verlag, New York.

\ref Stamey, T.A., Kabalin, J.N., McNeal, J.E., Johnstone, I.M., Freiha, F., Redwine, E.A. and Yang, N. (1989). Prostate specific antigen in the diagnosis and treatment of adenocarcinoma of the prostate: II. radical prostatectomy treated patients, Journal of Urology 141(5), 1076–1083.

\ref Tibshirani, R. (1996). Regression shrinkage and selection via the LASSO, {\em J. Royal. Statist. Soc. B.}, {\bf 58}(1), 267-288.

\ref Xu, H. K. (2014) Properties and iterative methods for the LASSO and its variants. {\em Chin. Ann. Math.} {\bf 35B}(3), 501--518.

\ref Zou, H. and Hastie, T. (2005). Regularization and variable selection via the elastic net.
{\em J. R. Stat. Soc. Ser. B Stat. Methodol.}, {\bf 67}(2), 301-–320.

%\ref Zou, H., Hastie, T., and Tibshirani, R. (2007). On the “degrees of freedom” of the lasso. {\em Ann. Statist.}, {\bf35}(35), 2173-–2192.
\end{document}